\title{Quasi-invariant means and Zimmer amenability}

\author{G\'abor Elek \and \'Ad\'am Tim\'ar\thanks{%
AMS Subject Classification: 43A07 \thinspace\ Research sponsored by  
Marie Curie Grant FP7 272955 and Sinergia grant CRSI22-130435 of the
Swiss National Science foundation}}

\documentclass[10pt,naturalnames]{article}
\usepackage{amsmath}
\usepackage{amsthm}
\usepackage{amsfonts}
\usepackage{graphicx}

\newtheorem{theorem}{Theorem}
\newtheorem{definition}{Definition}
\newtheorem{lemma}{Lemma}[section]

\newtheorem{proposition}{Proposition}[section]

\newtheorem{question}{Question}

\newcommand{\R}{\mathbb{R}}

\newcommand{\N}{\mathbb{N}}

\def \e {\epsilon}

\def \P {{\bf P}}

\def\cG{\mathcal{G}} \def\cS{\mathcal{S}}

\def \E {{\bf E}}

\def \proof {{ \medbreak \noindent {\bf Proof.} }}

\def \_reg {\rightarrow_{\bf reg}}

\def\maxdeg/{\Delta}

\def\wmu{\widetilde{\mu}}

\begin{document}
\maketitle

\bigskip

\begin {abstract}
Let $\Gamma$ be a countable group acting on a countable set $X$ by 
permutations. We give a necessary and sufficient condition for the action 
to have a quasi-invariant mean with a given cocycle. This can be viewed as
a combinatorial analogue of the condition for the existence of a quasi-invariant
measure in the Borel case given by Miller.
Then we show a geometric 
condition that guarantees that the corresponding action 
on the Stone-\v{C}ech compactification is Zimmer amenable. 
The geometric condition (weighted hyperfiniteness) resembles Property A. 
We do not know the exact relation between the two notions, however, we can 
show that amenable groups and groups of finite asymptotic dimension 
are weighted hyperfinite. 
\end {abstract}
\section{Introduction}
\noindent
{\bf Quasi-invariant means.}
Let $\Gamma$ be a countable group acting on a countable set $X$ by permutations.
Following John von Neumann we call a finitely 
additive probability measure on $X$ an {\it invariant mean} $\mu$ if 
it is preserved by 
the action, that is,
$$\mu(gA)=\mu(A)$$ holds for any $g\in\Gamma$ and $A\subseteq X$. 
The existence of the invariant mean is equivalent to the existence of 
a F\o lner sequence
$\{F_n\}^\infty_{n=1}$ having the following properties :
\begin{itemize}
\item
$F_n\subset \Gamma,\quad |F_n|$ is finite for any $n\geq 1$.
\item
For any $\epsilon>0$ and finite subset $L\subset\Gamma$ there 
exists a positive integer $n_{\e,L}$ such that if $n\geq n_{\e,L}$ then
$$\frac{|gF_n\cup F_n|}{|F_n|} < 1+\e $$
provided  $g\in L$.
\end{itemize}
Recall that a group is {\it amenable} if the natural left action on itself
possesses an invariant mean. Our first goal is 
to investigate quasi-invariant means, that is, finitely additive probability
 measures for which the zero measure class is preserved.

Before getting into details, let us consider the Borel analogue of the problem.
Let $\Gamma$ be a countable group acting by Borel automorphisms on a Polish
space $Y$. A Borel probability measure $\nu$ on $Y$ is an invariant measure if
$$\nu(gB)=\nu(B)$$
for any $g\in \Gamma$ and Borel set $B\subset Y$.
Miller studied quasi-invariant measures given by a given cocycle 
$\rho$ (see \cite{Miller}) for definitions). Recall that  a $\rho$-invariant
mean $\nu$ satisfies
$$\nu(gB)=\int_B \rho(g,x)\,d\nu(x)$$
for any $g\in \Gamma$ and Borel set $B\subseteq Y$. 
Miller gave a necessary and sufficient condition for
the existence of $\rho$-invariant means. 

Now let us consider the discrete analogue of $\rho$-invariant means.
Again, let $\Gamma$ be a countable group acting on a countable set $X$ by
permutations and $w$ be a positive real function on $X$. We further suppose that
for any fixed $g\in\Gamma$ the function $x\to \frac{w(gx)}{w(x)}$ is 
bounded on $X$ (later referred to as {\it boundedness condition}). First observe that $\rho(g,x):=\frac{w(gx)}{w(x)}$ 
is a $\Gamma$-cocycle on $X$. 
We say that a finitely additive probability measure $\mu$ on $X$ 
is $w$-invariant if for any $g\in\Gamma$ and $A\subseteq X$
\begin{equation} 
\label {i1}\mu(gA)=\int_A \frac{w(gx)}{w(x)} d\mu(x)\,\end{equation}
Recall \cite{Pat} that any finitely additive probability measure $\mu$ extends 
to a bounded linear functional on $l^\infty(X)$, hence the integral notation 
is meaningful. Since the step functions are dense in $l^\infty(X)$, 
(\ref{i1}) can be reformulated the following way:
\begin{equation} \label {i2}
\int_X F(gx)\,d\mu(x)=\int_X \frac{w(g^{-1}x)}{w(x)} F(x) d\mu(x)\,,
\end{equation}
where $F\in l^\infty(X), g\in\Gamma$. We will see that if $\Gamma$ is any 
finitely generated group then 
there exists $w:\Gamma\to\R^+$ satisfying the boundedness condition such that
$w$-invariant means exist on $\Gamma$ for the natural left action.
\begin{definition}
A family of finite subsets of $X$, $\{F_n\}^\infty_{n=1}$ forms a $w$-F\o lner
sequence if for any $\e>0$ and finite set $L\subset \Gamma$ there exists 
$n_{\e,L}$ such that if $n\geq n_{\e,L}$ and $g\in L$
then
$$\frac{\sum_{x\in gF_n\cup F_n} w(x)}{\sum_{x\in F_n}w(x)} < 1+\e\,.$$
\end{definition}
Notice that if $w=1$ then $w$-F\o lner sequences are exactly the usual
F\o lner sequences.
Our first theorem generalizes the classical result on the existence of invariant means.
\begin{theorem} \label{t1}
Let $\Gamma$ be a countable group acting on a countable set $X$ by permutations.
Let $w$ be a positive function on $X$ satisfying the boundedness condition.
Then the following two conditions are equivalent.
\begin{itemize}
\item There exist $w$-invariant means.
\item There exists a $w$-F\o lner sequence.

\end{itemize} \end{theorem}

\vskip 0.2in
\noindent
{\bf Weighted hyperfinite graphs.} Let $G$ be a connected, infinite graph
of bounded vertex degrees. We say that $G$ is {\it weighted hyperfinite} if
for any $\e>0$ there exists $K_\e>0$ with the following property.
For any finite induced subgraph $L\subset G$ and nonnegative function
$w:V(L)\to\R$ one can delete a subset $M$ of vertices 
(together with all the incident edges) of $L$ such
that
\begin{itemize}
\item $\sum_{x\in M} w(x)\leq \e \sum_{x\in V(L)} w(x)$
\item All the remaining components have size at most $K_\e$.
\end{itemize}
In this case we say that $M$ is a $(w,\epsilon)${\it -separating set} 
for $L$ with component sizes at most $K_\e$. If $H$ is a finite subset of
the vertices we define $w(H):=\sum_{x\in H} w(x)\,.$
We will call a positive function $w$ {\it balanced} if there exists
$C>1$ such that
$$\frac{1}{C}<\frac{w(y)}{w(x)}<C$$
for any adjacent vertices $x$, $y$. Observe that the boundedness condition
for $w$ is equivalent to being balanced. 
We will show that being weighted hyperfinite is invariant
under quasi-isometries. Hence, we can speak about weighted
hyperfinite groups. We will prove the following theorem.
\begin{theorem} \label{t2}
Let $\Gamma$ be a finitely generated amenable group or a group
of finite asymptotic dimension. Then $\Gamma$ is weighted hyperfinite.
\end{theorem}

\noindent
{\bf Zimmer amenability.}
Let $\Gamma$ be a finitely generated group with symmetric generating
system $S$. Let $\Gamma$ act 
 on a compact Hausdorff space $Y$ preserving
the measure class of a Borel probability measure. That is, the action
of $\Gamma$ is quasi-invariant. The orbit equivalence relation $E$ of
the action is defined the following way: $(xEy)$ is and only if $x$ and
$y$ are in the same orbit, that is, $gx=y$ for some $g\in\Gamma$.
The action is called {\it Zimmer amenable} \cite{EG},\cite{Zim} if
\begin{itemize}
\item It is hyperfinite, that is, there exist finite Borel subrelations 
(all the equivalence classes are finite)
$F_1\subset F_2\subset\dots$ such that $\cup^\infty_{n=1} F_n=E$ modulo 
a zero class.
\item Almost all the point stabilizers are amenable.
\end{itemize}
In this paper we will consider essentially free actions. In this case,
hyperfiniteness and Zimmer amenability coincide.

\noindent
Now let $\beta\Gamma$ be the Stone-\v{C}ech compactification of $\Gamma$. 
The elements of $\beta\Gamma$ are the ultrafilters on $\Gamma$. The principal
ultrafilters are identified with the elements of $\Gamma$. A base
of compact, Hausdorff topology on $\beta\Gamma$ is given by
$\{A_*\}_{A\subset\Gamma}$, where
$A_*$ is the set of ultrafilters containing $A$. Then
\begin{itemize}
\item $A_*\cup B_*=(A\cup B)_*$
\item $A_*\cap B_*=(A\cap B)_*$
\item $\overline{A}_*=\overline{(A_*)}$
\item $\emptyset_*=\emptyset$
\item $\Gamma_*=\beta\Gamma\,.$
\end{itemize}
Note that the sets $A_*$ are both closed and open.
Since $\Gamma$ acts on the ultrafilters of $\Gamma$, a continuous
$\Gamma$-action is given on $\beta\Gamma$, where
$$g(A_*)=(gA)_*$$
holds for any $g\in \Gamma$ and $A\subset\Gamma$.
Recall that the space of continuous functions on $\beta\Gamma$ can be identified
with $l^\infty(\Gamma)$. If $F\in l^\infty(\Gamma)$, we denote by
$F_*$ the corresponding element in $C(\beta\Gamma)$. Then
\begin{itemize}
\item $(FG)_*=F_* G_*$
\item $(F+G)_*=F_*+G_*$
\item $F_*\circ g=(F\circ g)_*$
\end{itemize}

By the Riesz representation theorem, if $\mu$ is a
finitely additive probability measure on $\Gamma$
then
$$\phi(F_*):=\int_\Gamma F d\mu$$
defines a regular Borel probability measure $\wmu$ on $\beta\Gamma$ such
that
$$\int_{\beta\Gamma} F_* d\wmu=\int_\Gamma F d\mu$$
for any bounded real function $F$ on $\Gamma$.

\noindent
Now let $w:\Gamma\to\R$ be a positive real function such that
for any $g\in\Gamma$ the real function on $\Gamma$ given by 
$x\to \frac{w(gx)}{w(x)}$ is bounded.
Define the function $z\to \rho(g,z)$ in $C(\beta\Gamma)$ as $(x\to \frac{w(gx)}{w(x)})_*$\,.
Then $\rho$ is a $\Gamma$-cocycle on $\beta\Gamma$. That is,
$$\rho(gh,z)=\rho(g,hz) \rho(h,z)\,.$$
Indeed, by the $\Gamma$-equivariance of the correspondance $F\to F_*$
\begin{itemize}
\item $z\to \rho(gh,z)=(x\to \frac{w(ghx)}{w(x)})_*$
\item $z\to \rho(g,hz)=(x\to \frac{w(ghx)}{w(hx)})_*$
\item $z\to \rho(h,x)=(x\to \frac{w(hx)}{w(x)})_*$
\end{itemize}

Therefore
$$\int_{\beta\Gamma} F_*(gz) d\wmu(z)= \int_{\beta\Gamma}
F_*(z)\rho(g^{-1},z) d\wmu(z)\,.$$
That is, $\wmu$ is a $\rho$-invariant measure. In other words, the action of
$\Gamma$ is quasi-invariant on $\beta\Gamma$ with Radon-Nykodym cocyle $\rho$.
The next result sheds some light on the relation between weighted 
hyperfiniteness and Zimmer amenability. 
 
\begin{theorem} \label{t3}
Let $\Gamma$ be a finitely generated group, with a positive, balanced
 weight function $w$
and a $w$-F\o lner sequence $\{F_n\}^\infty_{n=1}$. If $\Gamma$ is weighted 
hyperfinite, then the corresponding $\Gamma$ action on $\beta\Gamma$ with 
respect to the measure $\wmu_{F,\omega}$ is Zimmer amenable.
\end{theorem} 

\noindent
Let us recall that a group $\Gamma$ has Property A if and only if
its canonical action on the Stone-\v{C}ech compactification is topologically
amenable \cite{HR}. Also, by \cite[Corollary 3.3.8]{AR} if $\nu$ is
a quasi-invariant measure with respect to a free topologically amenable action
then it is Zimmer amenable. Hence, if $\Gamma$ is of Property $A$ then
the conclusion of Theorem \ref{t3} also holds. It is well-known that  
amenable groups as well as groups of finite asymptotic dimension
have Property A (by Theorem \ref{t2} they are weighted hyperfinite as well). 
Finally, if $\Gamma$ is finitely generated and
contains an embedded expander sequence, then it cannot have Property $A$
\cite{Roe}.
Clearly, $w$-hyperfinite groups cannot have imbedded expander sequences.
\begin{question}
What is the relation between weighted hyperfiniteness and Property A\,?
Does any of these properties imply the other\,?
\end{question}

\section{Quasi-invariant means}
Let $\Gamma$ be a finitely generated group
acting on the countable set $X$ with
a symmetric generating set $S$. Let $w$
be a positive, balanced real function on $X$. 
\begin{proposition} \label{p21}
If there exists a $w$-F\o lner sequence,
then there exist $w$-invariant means, as well.
\end{proposition}
\proof
Let $\{F_n\}^\infty_{n=1}$ be a $w$-F\o lner system and $\omega$ be a 
nonprincipal ultrafilter on $\N$. Let $\lim_\omega:l^\infty(X)\to \R$
be the corresponding ultralimit.
Define $\mu$ by
$$\mu(A):=\lim_\omega \frac{\sum_{A\cap F_n} w(x)}{\sum_{x\in F_n} w(x)}\,.$$
Then $\mu$ is clearly a finitely additive measure extending
to a mean on $l^\infty(X)$ by

$$\int_X F(x) d\mu(x):=
 \lim_\omega \frac{\sum_{x\in F_n} F(x)w(x)}{\sum_{x\in F_n}w(x)}\,.$$

Hence
$$\int_X F(gx) d\mu(x)=
 \lim_\omega \frac{\sum_{x\in F_n} F(gx)w(x)}{\sum_{x\in F_n}w(x)}\,.$$

By the $w$-F\o lner property
$$\lim_{n\to\infty} \frac{|\sum_{x\in F_n} F(gx)w(x)-
\sum_{x\in F_n} F(x)w(g^{-1}x)|}
{\sum_{x\in F_n}w(x)}=0$$
holds for any $g\in\Gamma$.

Hence,
$$\int_X F(gx) d\mu(x)=\lim_\omega\frac {\sum_{x\in F_n} F(x)w(g^{-1}x)}
{\sum_{x\in F_n}w(x)}= \int_X F(x) \frac{w(g^{-1}x)}{w(x)} d\mu(x)\,.$$
\qed

\vskip 0.2in
\noindent
We will denote the invariant mean constructed above by $\mu_{F,\omega}$.
Let $T\subset\Gamma$ be a finite set and
$\{\Psi_g\}_{g\in T}$ be bounded positive functions on $X$.
We say that $\{\Psi_g\}_{g\in T}$ is a $w$-compression system if for
any $x\in X$

$$\sum_{g\in T} \Psi_g(x)=1\quad\mbox{and}\quad \sum_{g\in T}
 \Psi_g (g^{-1}x) \frac{w(g^{-1}x)}{w(x)}<\frac{1}{2}\,.$$

The notion of a $w$-compression system is motivated by Miller's idea
of $\rho$-compressability \cite{Miller}.

\begin{proposition} \label{p22}
If there exists a $w$-compression system then there exists no
$w$-invariant mean.
\end{proposition}
\proof
Suppose that $\mu$ is a $w$-invariant mean. Then

$$1=\sum_{g\in T} \int_X \Psi_g(x) d\mu(x)
=\sum_{g\in T}  \int_X \Psi_g(g^{-1}x)\frac{w(g^{-1}x)}{w(x)}d\mu(x)
<\frac{1}{2}\,$$
providing a contradiction. \qed

\begin{proposition} \label{p23}
If $w$-F\o lner systems do not exist, then we have a $w$-comp\-ression system.
\end{proposition}
\proof 
Let $a_g(x):=\Psi_g(x)w(x)\,.$ Then for any $x\in X$
\begin{equation} \label{e8}
\sum_{g\in T} a_g(x)=w(x) \quad\mbox{and} \quad
\sum_{g\in T} a_g(g^{-1}x)<\frac{1}{2} w(x)\,.
\end{equation}

Thus we need to prove that there exists a system $\{a_g\}_{g\in T}$
satisfying (\ref{e8}).

A classical application of the Max Cut-Min Flow Theorem is the 
Transportation Problem. Say, we have a finite bipartite graph
$G=(V,E)$, where the vertex set $V$ is the disjoint union of $A$ and $B$.
Let $p:V\to\R$ be a positive function. One can think about the elements $x$ of
 set 
$A$ as manufacturer producing a certain (divisible) product worth of $p(x)$.
 The elements $y$ of $B$ are buyers having $p(y)$ amount of money to spend. 
A manufacturer can sell goods only to the buyers he is connected to. 
The question (the Transportation Problem) is whether the manufacturer can sell 
all the goods to the buyers or not. In mathematical terms the problem is
 to associate nonnegative numbers $a(x,y)$ to the edges $(x,y)$ such that
\begin{itemize}
\item
For any $x\in A$ 
$$\sum_{(x,y)\in E} a(x,y)=p(x)\,.$$
\item
$$\sum_{(x,y)\in E} a(x,y)\leq p(y)\,.$$
\end{itemize}

According to the Max Cut Min Flow Theorem, 
the sufficient and necessary condition
for the solvability of the transportation problem is that for any subset
$L\subseteq A$
$$\sum_{x\in L} p(x)\leq \sum_{y\in K} p(y)\,,$$
where $K$ is the set of vertices in $Y$ adjacent to a vertex in $X$.

\noindent
By compactness, the solvability of the Transportation Problem has the
same necessary and sufficient condition even if $G$ is an infinite bipartite
graph with bounded vertex degrees.
Let us see, how can we use the transportation problem.
Suppose that there exists no $w$-F\o lner system. Then there exists
a positive $\e>0$ such that for any finite subset $C\subset \Gamma$
$$(1+\e)\sum_{x\in C} w(x) < \sum_{y,\,y=gx, g\in S,x\in C}w(y)\,.$$
Hence there exists some $k>0$ such that for any $x\in X$
\begin{equation} \label{e10}
2\sum_{x\in C} w(x)< \sum_{y,\,y=gx, g\in S^k,x\in C}w(y)
\end{equation}
doubling condition is satisfied. Let us remark that the idea of using
doubling conditions  is due to Deuber, Simonovits and S\'os \cite{DSS}. 
In their paper they used the Marriage Lemma, which is also a classical
special case of the Max Flow Min Cut Theorem.

\noindent
Now let us construct our bipartite graph $G$. Let both the left and the 
right vertex set of $G$ be $X$. Draw an edge $(x,y)$ if
$y=g^{-1}x, g\in S^k:=T$. For the vertices $x$ on the left, define $p(x)$ to be
 $w(x)$. For the vertices $y$ on the right, define $p(y)$ to be 
$\frac{1}{2}w(y)$.
Then the equation (\ref{e10}) is just the necessary and sufficient condition
of the corresponding Transportation Problem. Hence (\ref{e8}) can be satisfied.
\qed

\vskip 0.2in
\noindent
By Propositions \ref{p21},\ref{p22} and \ref{p23}, Theorem \ref{t1} holds
if $\Gamma$ is finitely generated. Now let $\Gamma$ be an arbitrary countable
group acting on $X$ and let $w$ be a positive, balanced real function on $X$. 
Clearly, $\Gamma$ has a $w$-F\o lner
system if and only if all of its finitely generated subgroups possesses a
 $w$-F\o lner system.  Also, if there exists a $w$-invariant mean 
for each finitely generated subgroup then there exists 
a $w$-invariant mean for $\Gamma$ as well.

\noindent
Indeed, let $a_1, a_2,\dots $ be an enumeration of the elements of $\Gamma$ and
$K_n\subset M_\Gamma$ be the set of invariant means with respect to the group
generated by the set $\{a_1,a_2,\dots,a_n\}$, where $M_\Gamma$ is the compact 
Hausdorff space of all means on $\Gamma$ \cite{Pat}\,.

Since $K_1\supset K_2 \supset \dots $ is a sequence of nonempty closed sets,
there exists a $w$-invariant mean
$\mu\in\cap^\infty_{n=1} K_n\,.$
This finishes the proof of Theorem \ref{t1}\,.\qed

\section{Weighted hyperfinite graphs}
\begin{proposition}
Weighted hyperfiniteness is invariant under \\ quasi-isometries.
\end{proposition}
\proof
Let $G_1$ and $G_2$ be quasi-isometric graphs with a uniform bound $d$ on 
their vertex degrees, and suppose that $G_2$ is weighted hyperfinite. 
We need to show that $G_1$ is weighted hyperfinite as well.
Let $\iota : G_1\to G_2$ be a map that satisfies 
$c^{-1} d_{G_1}(x,y)- c \leq d_{G_2}(\iota (x),\iota (y)) \leq c d_{G_1}(x,y)+c$ 
with some $c>0$ for any $x,y \in G_1$,
and suppose that for every $z \in G_2$ there is some
$x \in G_1$ such that
$
d_{G_2}(z,\iota(x)) \leq c$ (note that slightly abusing notation we denote the 
graphs and their vertex sets by the same letters). In particular, 
for any $v\in V(G_2)$, $|\{x\,:\, \iota (x)=v\}|\leq c^2$. 
Note that $d^{2c+1}$ is an upper bound for the size of any ball of radius $2c$,
and fix $C:=\max \{d^{2c+1}, c^2\}$. 
Define a map $f: G_2\to G_1$ as $f(z):=x$, where $x$ is a point that 
minimizes $d_{G_2}(z,\iota (x)$
(fixed arbitrarily, in case of ambiguity). In particular, $\iota (f(z))=z$ 
when $z\in\iota (G_1)$. We mention that every point of $G_1$ has 
at most $C^2$ preimages by $f$.

Let $H_1$ be an arbitrary finite induced subgraph of $G_1$; we want to 
show that for any $\epsilon >0$ and any weight function $w$ on $H_1$, 
there is a $(w,\e)$-separating set for $H_1$ with 
component sizes independent of the choice of $H_1$. 
Define $H_2=\iota (H_1)$, and let $H_2^+$ be 
the $2c$-neighborhood of $H_2$. Define a weight function $w'$ 
on $H_2^+$ to be $w'(z):=w(f(z))$. Then we have
\begin{equation}
\label{elso}
w'(H_2^+)\leq C^2 w(H_1)
\end{equation}
by our observation on the number of preimages by $f$.
Now define $w''$ on $H_2^+$ 
by letting $w''(z):=\sum_{y\in B(z)} w'(y)$, where $B(z)$ is the $2c$-neighborhood of $z$ in $H_2^+$. We have noted that $|B(z)|\leq C$, hence
\begin{equation}
\label{masodik}
w''(H_2^+)\leq Cw'(H_2^+)\leq C^3 w(H_1),
\end{equation}
using (\ref{elso}) for the second inequality.

Let $S$ be an $(w'',\e)$-separating set for $H_2^+$ with components 
of sizes $K(\epsilon )$, and let $S^+$ be the $2c$-neighborhood of $S$ in $H_2^+$. By definition of $w''$ we have
\begin{equation}
\label{harmadik}
w'(S^+)\leq w''(S).
\end{equation}

We claim that $\iota^{-1} (S^+) $ is an $(w,C^3\epsilon)$-separating set for $H_1$ of component sizes $CK(\epsilon)$ (which would complete the proof, since $\epsilon$ was arbitrary, and $C$ only depended on $c$ and $d$). First, $w (\iota^{-1} (S^+))\leq w' (S^+)\leq w''(S)$ using (\ref{harmadik}), and $w(H_1)\geq C^{-3} w'' (H_2^+)$ by (\ref{masodik}). Thus $w (\iota^{-1} (S^+))/w(H_1)
\leq C^3 w''(S)/w''(H_2^+)\leq C^3\epsilon$. 
So it only remains to show that the components of $H_1\setminus \iota^{-1}(S^+)$ have sizes at most $CK(\epsilon)$.
This follows from the next claim.

{\it Claim:} Let $\iota (x),\iota (y)\in H_2$ be in different components of $H_2^+\setminus S$. Then $x$ and $y$ are in different components of $H_1\setminus \iota^{-1} (S^+)$.

Suppose not, and let $P$ be a path between $x$ and $y$ in $H_1\setminus \iota^{-1} (S^+)$. Consider $\iota (P)$. Since $P\cap \iota^{-1} (S^+)=\emptyset$, $\iota (P)\cap S^+=\emptyset$. Hence $\iota (P)$ is at distance at least $2c$ 
from any element of $S$. On the other hand, two consecutive (adjacent) vertices $u$ and $v$ in $P$
are mapped into points at distance at most $2c$ by the quasi-isometry $\iota$. Therefore we can connect each such pair $\iota (u), \iota (v)\in H_2$ by a path of length at most $2c$ in $H_2^+$, which path is thus disjoint from $S$. The union of these paths between $\iota(u),\iota(v)$ in $H_2^+$ (over all such $u,v$) avoids $S$, hence $\iota (x)$ and $\iota(y)$ are in the same component of $H_2^+\setminus S$. This contradicts the assumption on $x$ and $y$, finishing the proof.
\qed

\vskip 0.2in
\noindent
{\bf Proof of Theorem \ref{t2}}.
Now let us prove that bounded degree graphs of finite asymptotic dimension are
weighted hyperfinite. Recall that a graph $G$ has asymptotic dimension $d$ if
for every $r>0$ there exists an $R(r)=R$ and vertex-disjoint induced subgraphs 
${\cal U}_1,\ldots , {\cal U}_d$ of $G$, such that 
every vertex of $G$ is in some ${\cal U}_i$, 
for each $i\in\{1,\ldots ,d\}$ every connected 
component of ${\cal U}_i$ has diameter at most $R$, and 
any two distinct components of ${\cal U}_i$ have distance at least $r$ 
in $G$. If $d$ is finite then we say that $G$ 
has {\it finite asymptotic dimension}. 
The asymptotic dimension is a quasi-isometry invariant 
(hence it defines a group invariant). See e.g. \cite{BD} for a 
survey on the asymptotic dimension.

\noindent
So, let $d$ be the asymptotic dimension of the bounded degree 
graph $G$. Let $d$ be the asymptotic dimension of $G$, $H$ be 
an arbitrary induced subgraph of $G$, $w: V(H)\to \R$ be a 
weight function on the 
vertices,
 and $\epsilon >0$. Define $r:= 2[1+1/\epsilon]$, and 
let ${\cal U}_1,\ldots , {\cal U}_d$ be the families of 
sets corresponding to $r$  in the definition of asymptotic 
dimension, and $R$ be the corresponding $R(r)$. 
For $i\in\{1,\ldots , d\}$, $t\in\{1, 2,\ldots , 1+[1/\epsilon]\}$, 
let $S_i(t)\subset V(G)$ be the set of 
points at distance $t$ from ${\cal U}_i$. 
In particular, the sets $S_i(1),\ldots S_i (1+[1/\epsilon])$ 
are pairwise disjoint. Hence 
there is a $t\in\{1, 2,\ldots , 1+[1/\epsilon]\}$ such 
that $w(S_i(t))\leq w(H)/(1+[1/\epsilon])\leq \epsilon w(H)$; 
let $t(i)$ be one such $t$. 
On the other hand, any two components of ${\cal U}_i$ 
are at distance at least $r= 2(1+[1/\epsilon]) $ 
from each other, thus any two such components are 
separated by $S_i(j)$ for any $j$. 
We obtain that $S:=\cup_{i=1}^d S_i (t(i))$ 
is such a set that any component of $H\setminus S$ 
intersects at most one component of each $U_i$, 
hence its total diameter is at most $d(R+1)$. 
The uniform bound on the degrees of $G$ then implies that 
every component of $H\setminus S$ has a 
uniformly bounded size. 
Finally, we have $w(S)= \sum_{i=1}^d w(S_i(t(i)))\leq d\epsilon w(H)$. 
Since $\epsilon$ was arbitrary, this shows that $G$ is indeed hyperfinite.

\vskip 0.2in
\noindent
Now we prove that the Cayley graph of a finitely generated amenable
group is weighted hyperfinite. Let $\Gamma$ be a finitely generated amenable
group with symmetric generating set $S$ and $G$ be its left-Cayley graph. That 
is the vertex set of $G$ is $\Gamma$ and the vertices $x$ and $y$ are connected
if $x=sy$, for some $s\in S$. Let $\{F_n\}^\infty_{n=1}$ be a F\o lner sequence 
in $\Gamma$.
For later convenience, suppose that $|\partial F_n|/|F_n|\leq n^{-2}$. Recall 
that $\partial F_n$ is the set of elements in $F_n$ that are connected to a 
vertex in the complement of $F_n$. Thinking about vertices of $G$ as elements of the group, we will refer to products of vertices and vertex sets. Also, assume that the identity element is contained in each of the $F_n$.
Note that for a subgraph $H$ and $g\in \Gamma$,
the map $x\mapsto xg$ from $H$ preserves edges.

Let $H$ be an arbitrary induced subgraph of $G$, and $w$ a weight function on
 its vertices. 
Set $p_n\in [0,1]$ to be such that $(1-p_n)^{|\partial F_n|}
=1-n^{-1}$. 
 
Define a random set $R_n$ in $G$ as follows: an $x\in V(G )$ will be in $R_n$ 
with probability $p_n$ and independently from the others. 
Define $B_n$ to be $B_n=V(H)\setminus F_nR_n$. Finally, let $S_n\subset H$ 
be defined as $(\partial F_n )R_n\cap V(H)$. 
Any component of $H\setminus (B_n\cup S_n)$ has size at most $|F_n|$, 
because if $x\in V(H)$ is not in $B_n$, then there is some $v\in R_n$ with $v\in F_n^{-1}x$, and hence $S_n$ separates $x$ from any point in
 the complement of $F_nv$. 

We claim that $B_n\cup S_n$ has relatively small weight 
for some choice of $R_n$.
First, its expected weight is:
$$\E[ w(S_n)+w(B_n)]=\E [\sum_{x\in H} w(x)  {\bf 1}_{x\in \partial F_n R_n}]
+\E [\sum_{x\in H} w(x)  {\bf 1}_{x\not\in F_n R_n}]=$$
$$= \sum_{x\in H} w(x) (\P[x\in \partial F_nR_n]+\P[x\not\in F_nR_n])
=$$
$$= \sum_{x\in H} w(x) (\P[\partial F_n^{-1}x\cap R_n\not=\emptyset]+
\P[F_n^{-1}x\cap R_n=\emptyset])=
$$
$$= w(H) (1-(1-p)^{|\partial F_n|}+(1-p)^{|F_n|})\leq w(H) (n^{-1}+e^{-n}),$$
where the last inequality follows from the assumption on $F_n$ and the 
choice of $p$. Hence, there is some $R_n$ where the corresponding $S_n,B_n$ 
satisfies $w(S_n\cup B_n)\leq 2n^{-1}w(H)$. 
We have also observed that $S_n\cup B_n$ splits $H$ into pieces of sizes 
at most $|F_n|$. Since $H$ was arbitrary, we have proved 
that $G$ is hyperfinite.
\qed

\section{The Proof of Theorem \ref{t3}}
In this section we use some ideas from \cite{AE}.
Let $\Gamma$ be a finitely generated group
with symmetric generating system $S$ and let $G$ be the associated
(left) Cayley graph. The graphing $\cG$ (see \cite{KM}) of the
associated $\Gamma$-action on $\beta\Gamma$ is defined as 
follows: $x,y\in\beta\Gamma$ are connected if there exists $s\in S$
such that $sx=y$.

\begin{lemma}
The action of $\Gamma$ on $\beta\Gamma$ is free. Hence the components
of $\cG$ are isomorphic to $G$.
\end{lemma}
\proof Let $g\in \Gamma$ and $\omega\in\beta\Gamma$ be an ultrafilter.
Then $g\omega$ is the ultrafilter containing the sets $gA$, where
$A\in\omega$. Let $\cup_{i=1}^n A^i=\Gamma$ be a finite partition
such that $gA^i\cap A^i=\emptyset$ for a fixed $g\in\Gamma$ 
and $1\leq i \leq n$.
Then $g(A^i_*)\cap A^i_*=\emptyset\,.$ Since $\beta\Gamma=\cup_{i=1}^n A^i_*$
$g$ cannot fix any element of $\beta\Gamma$. \qed

\vskip 0.1in \noindent
Let $T\subseteq G$ be a subgraph such that the vertex set of $T$ is 
the whole $\Gamma$. We can associate a Borel subgraphing
$\cG(T)\subset\cG$ to $T$ the following way.
For $s\in S$, let $A_s\subset\Gamma$ be the set of vertices $x$ 
such that $x$ and $sx$ are adjacent in $T$. Now connect $y\in (A_s)_*$ to
$sy$. Hence the subgraphing $\cG(T)$ is the union of the graphs
of the Borel automorphisms $(A_s)_*\to s(A_s)_*$, where $s\in S$.
Note that by a Borel subgraphing (as in \cite{KM}) we always mean a
Borel subgraph of the graphing, such that the vertex set is the whole
space $Y$.
\noindent
If $T\subset G$ let $\partial T$ be the set of vertices $x$ such that
$x$ is adjacent to some $y$, that is, not in the same $T$-component as
$x$. Similarly, we can define $\partial(\cG(T))\,.$

\begin{proposition}
Let $\Gamma,w$ and $\{F_n\}^\infty_{n=1}$ be as in Theorem \ref{t3}. Then
for any $\e>0$ there exists a subgraph $T_\e\subset G$ with components
of bounded size such that 

\begin{equation} \label{e51}
\mu_{F,\omega}(\partial T_\e)<\e\,.
\end{equation} \end{proposition}
\proof
Since $\{F_n\}^\infty_{n=1}$ is $w$-F\o lner, for any $\delta>0$ we have
a subgraph $F^\delta_n$, $\frac{w(F^\delta_n)}{w(F_n)}<\delta\,,$
such that if we delete $F^\delta_n$ from $F_n$ the components of the resulting
graph are bounded by $K_\delta$. Let $T_\delta$ be the union
of these components plus all the points outside the union of 
the $w$-F\o lner sets as singletons. 
That is the vertex set of $T_\delta$ is $\Gamma$.  Clearly,
if $x\in\partial T_\delta\cap F_n$ then either
$x\in\partial F_n$ of $x$ is in the $1$-neighborhood of $F^\delta_n$,
$B_1(F^\delta_n)$.
By the F\o lner property, 
$$\lim_{n\to\infty} \frac{w(\partial F_n)}{w(F_n)}=0\,.$$
Also,
$$w(B_1(F^\delta_n))\leq C|S| w(\partial F_n)\,,$$
where $C$ is the constant in the balancedness condition for $w$.
Therefore,\\ $\mu_{F,\omega}(\partial T_\e)\leq C|S|\delta\,.$
Hence if $\delta<\epsilon/C|S|$ the equation (\ref{e51}) is satisfied. \qed

\begin{lemma}
Let $T\subset G$ be a subgraph
such that all the components of $T$ have  size at most $k$. Then all
the components of $\cG(T)$ have size at most $k$ as well. Moreover,
$\partial (\cG(T))=(\partial T)_*\,.$ \end{lemma}
\proof
Let $A\subseteq \Gamma$ be a set containing exactly one element
from each component. For each subgraph of $G$ we have a set 
$A_F\subset A$ of vertices $x$ such that
the component of $x$ is
isomorpic to $Fa$ (even as a graph with edge labels from $S$).
Then 
\begin{equation} \label{e141}
\cG(T)=\cup_F F(A_F)_*\,,
\end{equation}
hence all the components of $\cG(T)$ have size at most $k$. For the
second statement,
$$\partial(\cG(T))= \partial (\cup_F F(A_F)_*)=
\cup_F \partial(F(A_F)_*)=\cup_F(\partial F)(A_F)_*=(\partial T)_*$$
\qed

\vskip 0.2in

By \cite[Proposition 10.3]{KM}, $\cG$ is hyperfinite 
if and only if for
any $\e>0$ there exists a Borel subgraphing $\cS_\e\subset\cG$ such that
\begin{equation}
\label{e661}
\wmu_{F,\omega}(\partial \cS_\e)\leq \e
\end{equation}
Let $\cS_e:=(T_\e)_*$, then (\ref{e661}) follows. This ends the
proof of Theorem \ref{t3}\,. \qed

\vskip 0.2in
\noindent
Finally, let us show that for any finitely generated group $\Gamma$,
 there exists a positive, balanced function $w:\Gamma\to\R$ such that:
\begin{itemize}
\item there exist $w$-F\o lner systems,
\item the resulting measure $\wmu_{F,\omega}$ is atomless.
\end{itemize}
First of all, we can suppose that $\Gamma$ is nonamenable, since for amenable
groups $w:=1$ clearly satisfies the two conditions.
Let $\{B_r(x_r)\}^\infty_{r=1} $ be vertex disjoint 
balls in a Cayley-graph $G$ of $\Gamma$. Note that the distance
$d$ in $G$ is the shortest path metric and $B_r(x_r)$ is the $r$-ball
around $x_r$. Let $w_r$ be defined on $B_r(x_r)$ the following way.
For any $0\leq i \leq r$, $ w_r(S_i(x_r))=\frac{1}{r+1}\,,$ and
$w_r(x)=w_r(y)$ if $x,y\in S_i(x_r)$, where 
$$S_i(x_r)=\{y\,\mid \,d(x_r,y)=i\}\,.$$
Then, for $x\in B_r(x_r)$ let $w(x)=\frac{w_r(x)}{w_r(z)}$,
where $z\in S_r(x_r)$. If $x\notin \cup^\infty_{r=1}B_r(x_r)$, let $w(x)=1$.
Then $w$ is balanced, since by nonamenability, and
by the fact that $\Gamma$ is finitely generated, there exists $D>1$ such
that for all $r\geq 1$ and $0\leq i \leq r$, 
$$\frac{1}{D}<\frac{|S_{i+1}(x_r)|}{|S_i(x_r)|}< D $$
Clearly, $\{B_r(x_r)\}^\infty_{r=1} $ forms a $w$-F\o lner system.
Observe that for any $k\geq 1$ there exists $r_k\geq 1$ such that
if $r\geq r_k$ then we can partition $B_r(x_r)$ into $k$ parts such that
the weight of each part is less than $\frac{2}{k} w(B_r(x_r))$.
This observation easily follows from the fact that
$$\lim_{r\to\infty} \frac{\max_{y\in B_r(x_r)}w(y)}{w(B_r(x_r))}=0\,.$$
Therefore, for any $\delta>0$ one can partition $\beta\Gamma$ into
finitely many Borel parts such that the $\wmu_{F,\omega}$-measure of each
part is less than $\delta$. Hence $\wmu_{F,\omega}$ is atomless.

\noindent
\'Ad\'am Tim\'ar, Fakult\"at f\"ur Mathematik, Universit\"at Wien
Nordbergstra\ss e 15, 1090 Wien
\vskip 0.2in
\noindent
G\'abor Elek, Alfred Renyi Institute, Realtanoda u. 13-15, 1053 Budapest and
EPFL Station 8
CH–1015 Lausanne

\end{document}